\documentclass[11pt]{article}
\usepackage{fontenc}
\usepackage{inputenc}
\usepackage{authblk}
\usepackage{mathtools}

\mathtoolsset{showonlyrefs}
\usepackage{amsmath}
\usepackage{amssymb}
\usepackage{amsthm}
\usepackage{subfigure}
\usepackage{caption2}
\usepackage{array}
\usepackage{xy}
\usepackage[pdftex]{graphicx}
\usepackage{hyperref}
\usepackage{color}
\usepackage{transparent}
\usepackage{latexsym}
\usepackage{float}
\usepackage{mathrsfs}
\usepackage{amssymb,url}

\makeatletter
\@addtoreset{equation}{section}
\makeatother

\newtheorem{lem}{Lemma}[section]
\newtheorem{thm}{Theorem}[section]
           
\newtheorem{cor}{Corollary}[section]
\newtheorem{definition}{Definition}[section]

\newtheorem*{rmk}{Remark}
\newtheorem*{pf}{Proof}

\usepackage{hyperref}

\usepackage{xcolor}
\begin{document}
	\title{\Large \bf The $L_p$ chord Minkowski problem for negative $p$}

	\author{ \small \bf Yuanyuan Li
		\\ \small  School of Mathematical Sciences, University of Science and Technology of China,\\
		\small Hefei, 230026, China \\  \small E-mails:
		lyuanyuan@mail.ustc.edu.cn
	}
    \date{}
	\maketitle
	
	\begin{abstract}
	In this paper, we solve the $L_p$ chord Minkowski problem in the case of discrete measures whose supports are in general position for negative $p$ and $q>0.$ As for general Borel measure with a density, we also give a proof but need $p\in(-n,0)$ and $n+1>q\geqslant 1.$ The $L_p$ chord Minkowski problem was recently posed by Lutwak, Xi, Yang and Zhang, which seeks to determine the necessary and sufficient conditions for a given finite Borel measure  such that it is the $L_p$ chord measure of a convex body, and it includes the chord Minkowski problem and the $L_p$ Minkowski problem.
	\end{abstract}

	\section{introduction}
	The central objects in study of convex geometry are convex bodies. A convex body in $n$-dimensional Euclidean space $\mathbb{R}^n,$ is a compact convex set with non-empty interior. The Brunn-Minkowski theory is a study of convex bodies which centers around the study of geometric functionals and the differential of these functionals. When geometric invariants arise as geometric functionals of convex bodies, geometric measures are often viewed as differentials of geometric invariants. One of the cornerstones of the Brunn-Minkowski theory is the Minkowski problem. It is a problem of priscribing geometric measure generated by convex bodies, which is concerned about necessary and sufficient conditions for a given measure such that it arises as the measure generated by a convex body. The most studied Minkowski-type problem is the classical Minkowski problem, which focuses on the surface area measures of convex bodies. For a comprehensive discussion on the Minkowski problem and its resolution, we recommend readers consulting  Pogorelov \cite{Pogorelov} and Cheng–Yau \cite{Cheng-Yau}.
	
	% When solution has good regularity, it's equivalent to solve a degenerate fully non-linear partial differential equation.
	Recently, a new family of geometric measures were introduced by Lutwak-XYZ\cite{LXZY2022} by studying of a variational formula regarding intergral geometric invariants of convex bodies called chord integrals. Let $K \in \mathcal{K}^n$ where $\mathcal{K}^n:=\{\text{all convex bodies in } \mathbb{R}^n\},$ the $q$th chord integral $I_q(K)$ is defined by
	\begin{equation}\label{integral}
	I_q(K)=\int_{\mathcal{L}^n}|K\cap\ell|^q d\ell,
	\end{equation}
	where $\mathcal{L}^n$ denotes the Grassmannian of 1-dimensional affine subspace of $\mathbb{R}^n,$ $|K\cap \ell|$ denotes the length of the chord $K\cap \ell,$ and the integration is with respect to Haar measure on the affine Grassmannian $\mathcal{L}^n,$ which is normalized to be a probability measure when restricted to rotations and to be $(n-1)$-dimensional Lebesgue measure when restricted to parallel translations.
	$$
	I_1(K)=V(k),\quad  I_0(K)=\frac{\omega_{n-1}}{n\omega_n}S(K), \quad I_{n+1}(K)=\frac{n+1}{\omega_n}V(K)^2,
	$$
	where $\omega_n$ denotes the volume of $n$-dimensional unit ball. Note that $I_q(B_n)=\frac{2^q\omega_n\omega_{n+q-1}}{\omega_{q}},$ where $B_n$ is the n-dimensional unit ball.  One can see from the above fomula that the chord integrals include the convex body's volume and surface area as two special cases. These are Crofton’s volume formula, Cauchy’s integral formula for surface area, and the Poincar\'{e}-Hadwiger formula,
	respectively (see [\cite{D. Ren}, \cite{Santalo}]).
	
	The chord measures and the Minkowski problems associated with chord measures were posed in \cite{LXZY2022}. They showed that the chord measures are the differentials of chord integrals and completely solved the chord Minkowski problem except for the critical case of the Christoffel-Minkowski problem. The $q$th chord measure is a finite Borel measure on $\mathbb{S}^{n-1}$ defined by
	\begin{equation}\label{measure1}
	F_q(K,\eta)=\frac{2q}{\omega_n}\int_{\nu^{-1}_K(\eta)}\tilde{V}_{q-1}(K,z)d \mathcal{H}^{n-1}(z),\text{   Borel }\eta\subset\mathbb{S}^{n-1},
	\end{equation}
	where $\widetilde{V}_{q-1}(K,z)$ is the $q-1$ th dual quermassintegral with respect to $z$.(See \eqref{dual}.)
	$$
	F_0(K,\cdot)=\frac{(n-1)\omega_{n-1}}{n\omega_n}S_{n-2}(K,\cdot), \quad F_1(K,\cdot)=S_{n-1}(K,\cdot),
	$$
	where $S_{i}(K,\cdot)$ is the $i$th order area measure of $K.$ Once chord measures are constructed, the $L_p$ chord measures follow naturally by extensions.
	For $K\in \mathcal{K}^n_{o}$ and $p\in \mathbb{R},$ the $L_p$ chord measures are defined by
	\begin{equation}\label{measure2}
	F_{p,q}(K,\eta)=\frac{2q}{\omega_n}\int_{\nu^{-1}_K(\eta)}(z\cdot\nu_{K}(z))^{1-p}\tilde{V}_{q-1}(K,z)d \mathcal{H}^{n-1}(z),\text{   Borel }\eta\subset\mathbb{S}^{n-1}.
	\end{equation}
	When $p=0,$ it is the cone-chord measure. When $q=1,$ $F_{p,1}(K,\cdot)$ is the $L_p$ surface area measure. When $q=0,$ $F_{p,0}(K,\cdot)$ is the $L_p-(n-2)$th area measure.
	
	The $L_p$-Minkowski problem was first formulated and studied by Lutwak in \cite{L1}. It has been rapidly attracting much attention; Lutwak introduced the important $L_p$ surface area measure and its associated Minkowski problem in the $L_p$ Brunn-Minkowski theory. Many cases of the $L_p$ Minkowski problem have been solved. The logarithmic Minkowski problem is one of the most central Minkowski type problems and is the problem of characterizing the cone-volume measure; see B$\ddot{o}$r$\ddot{o}$czky, Lutwak, Yang and Zhang [\cite{BLYZ13}, \cite{BLYZ12(5)}], Zhu [\cite{Zhu14},\cite{BHZ16}], Stancu [\cite{Stancu0240}, \cite{Stancu41}], Gage \cite{Gage18}, Xi and Leng \cite{XL43}, Firey \cite{Firey(17)}, Andrews \cite{Andrews (1)}, Chen, Huang, Li and Liu \cite{CHLL(10)}, Chen, Feng, Liu \cite{CFL 2022}, [\cite{H.Y45}, \cite{BCD(7)}] and reference therein. The centro-affine Minkowski problem is unsolved, see \cite{WZ06}. For more classical Brunn-Minkowski theory and its recent developments, we suggest readers to Schneider's book \cite{Schneider}.
	
	The $L_p$ chord Minkowski problem posed by Xi-LZY \cite{LXZY2022} is a problem of prescribing the $L_p$ chord measure: Given a finite Borel measure $\mu$ on $\mathbb{S}^{n-1},p\in \mathbb{R},$ and $q\geqslant 0.$ Asking what are the necessary and sufficient conditions for $\mu$ such that $\mu$ is the $L_p$ chord measure of a convex body $K\in \mathcal{K}^n_o,$ namely
	\begin{equation}\label{original eq}
	F_{p,q}(K,\cdot)=\mu
	\end{equation}
	when $p=1,$ it is the chord Minkowski problem. When $q=1,$ it is the $L_p$ Minkowski problem.
	When $\mu$ has a density $f$ that is an integrable nonnegtive function on $\mathbb{S}^{n-1},$ equation\eqref{original eq} becomes a new type of Monge-Amp\`{e}re equation on $\mathbb{S}^{n-1}$:
	\begin{equation}\label{MAeq}
	\mbox{det}(\nabla^2h+hI) =\frac{h^{p-1}f}{\widetilde{V}_{q-1}([h],\bar{\nabla} h)},\text{  on }\mathbb{S}^{n-1},
	\end{equation}
	where $\nabla^2h$ is the covariant differentiation of $h$ with respect to an orthonormal frame on $\mathbb{S}^{n-1},$ we look for a solution $h$ which is the support function for some nondegenerate convex body. We can extend $h$ to $\mathbb{R}^n$ via homogeneity and $\bar{\nabla} h$ is the Euclidean gradient of $h$ in $\mathbb{R}^n,$ and $\widetilde{V}_{q-1}([h],\bar{\nabla} h)$ is the $(q-1)$th dual quermassintegral of the Wulff-shape $[h]$ of $h$ with repect to the point $\bar{\nabla} h.$ 
	
	In their paper\cite{LXZY2022}, Lutwak-XYZ gave a sufficient condition for the symmetric case of the chord log-Minkowski problem by studying the delicate concentration properties of cone-chord measures. Shortly thereafter, Xi, Yang, Zhang and Zhao \cite{XYZZ2022} solved the $L_p$ chord Minkowski problem for $ p>1$ and for $0<p<1$ under the symmetric condition, where the origin symmetry played a crucial role in the case of $0< p< 1.$ More recently, Xi, Guo and Zhao solved the $L_p$ chord Minkowski problem when $0 \leqslant p < 1$, without any symmetry assumptions.
	
	In this paper, we say a polytope is a convex hull of finitely many points. The approximation of general convex bodies by simpler ones such as polytopes or bodies with differentiable boundaries is a useful tool for many investigations. Hence, polytopes are of great importance for the Minkowski problem and the $L_p$ Minkowski problem.[see, Lutwak-YZH\cite{LYZH2005}; Schneider's book\cite{Schneider} 392-393] By $\mathcal{P}$ we denote the set of all polytopes in $\mathbb{R}^n$, and a i-dimensional face of $P \in \mathcal{P}$ is referred to as an i-face. We call a face of dimension dim $P$-1 is a facet.
	\begin{definition}
		A finite set $U$ of unite vectors in $\mathbb{R}^n$ is said to be in general position if $U$ is not contained in a closed hemisphere of $\mathbb{S}^{n-1}$ and any $n$ elements of $U$ are linearly independent.
	\end{definition}
	Károlyi and Lovász \cite{KL polytope} were pioneers in the study of polytopes with facet normals in general position. Such polytopes are of great significance as any convex body can be approximated by a sequence of polytopes with facet normals in general position.
	Now we state our main theorems.
	\begin{thm}\label{A}
		Let $p<0,q>0,$ and $\mu$ be a discrete measure on $\mathbb{S}^{n-1}$ whose outer unit normals are in general position in dimension $n.$ Then there exists a polytope $P$ containing the origin in its interior such that
		\begin{equation*}
		F_{p,q}(K,\cdot)=\mu.
		\end{equation*}
	\end{thm}
    \begin{thm}\label{B}
    	Let $p\in (-n,0),n+1>q\geqslant 1,$ and $\mu$ be a finite Borel measure on $\mathbb{S}^{n-1}$ with density $f\in L^{\infty}(\mathbb{S}^{n-1}),$ and $\frac{1}{\Lambda}<f<\Lambda$ for some constant $\Lambda>0,$ then there exists a convex body $K\in \mathcal{K}^n_o$ such that 
    	\begin{equation*}
    	F_{p,q}(K,\cdot)=\mu.
    	\end{equation*}
    \end{thm}
    
    The remainder of this paper is structured as follows. In Section 2, we present fundamental concepts in the theory of convex bodies and integral geometry. Section 3 and 4 are dedicated to proving Theorems \ref{A} and \ref{B}, respectively.
    \section{Preliminaries}
    In this section, our objective is to establish notations and gather relevant results from the literature that will be necessary for the subsequent analysis. 
    
    We denote $x\cdot y$ as the standard inner product of $x,y\in \mathbb{R}^n,$ and write $|x|=\sqrt{x\cdot x}$ for the Euclidean norm of $x.$ We write $\mathbb{S}^{n-1}$ as $(n-1)$-dimension unit sphere of $\mathbb{R}^n,$ and denote $\mathcal{H}^{n-1}$ as the $(n-1)$-dimensional spherical Lebesgue measure. Denote $\mathcal{K}^n$ for the collection of all convex bodies in $\mathbb{R}^n$ and $\mathcal{K}^n_o$ for the subset of $\mathcal{K}^n$ that contains the origin in the interior.
    
    Let $\Omega\subset \mathbb{S}^{n-1}$ be a closed set of the unit sphere, not lying in a closed hemisphere, and a positive continuous function $h: \mathbb{S}^{n-1}\rightarrow \mathbb{R}$ is given.(Only the values of h on $\Omega$ will be needed, but without loss of generality we may assume that $h$ is defined on all of $\mathbb{S}^{n-1}.$) The Wulff shape of $h$ is defined by
    $$
    [h]=\{x\in\mathbb{R}^n: x\cdot u\leqslant h(u)\text{ for all } u\in\mathbb{S}^{n-1}\}.
    $$

    Let $K\in \mathcal{K}^n,$ $h(v)=h_K(v)=\max\{v\cdot x,x\in K\},$ $\rho(u)=\rho_K(u)=\max\{\lambda:\lambda u\in K\}$ are the support function and the radial function of convex body $K$ defined from $ \mathbb{S}^{n-1}\rightarrow \mathbb{R}.$ We write the support hyperplane of $K$ with the outer unit normal $v$ as
    $$
    H_K(v)=\left\{x\in \mathbb{R}^n: x\cdot v=h(v) \right\},
    $$
    the half-space $H^{-}(K,v)$ in direction $v$ is defined by
    $$
    H^{-}_K(v)=\left\{x\in \mathbb{R}^n: x\cdot v\leqslant h(v) \right\}.
    $$  
    Denote $\partial K$ as the boundary of $K$, that is, $\partial K=\{\rho_K(u)u:u\in \mathbb{S}^{n-1}\}.$ The spherical image $\nu=\nu_{K}:\partial K\rightarrow \mathbb{S}^{n-1}$ is given by
    $$
    \nu(x)=\{v\in\mathbb{S}^{n-1}:x\in H_K(v)\},
    $$
    let $\sigma_K\subset \partial K$ denote the set of all points $x\in \partial K,$ such that the set $\nu_K(x)$ contains more than one element. Fortunately, we have $\mathcal{H}^{n-1}(\sigma_K)=0$ (see \cite[page 84]{Schneider} ) 
    and the radial Gauss image $\alpha=\alpha_K$ and the reverse radial Gauss image $\alpha^*=\alpha^*_K$ are respectively defined by
    $$
    \alpha(\omega)=\{\nu(\rho_K(u)u):u\in \omega\},\alpha^*(\omega)=\{u\in \mathbb{S}^{n-1} \nu(\rho_K(u)u)\in \omega\}.
    $$
    Let $K \in \mathcal{K}^n$, for $z \in \operatorname{int} K$ and $q \in \mathbb{R}$, the $q$ th dual quermassintegral $\widetilde{V}_q(K, z)$ of $K$ with respect to $z$ is defined by
    \begin{equation}\label{dual}
    \widetilde{V}_q(K, z)=\frac{1}{n} \int_{S^{n-1}} \rho_{K, z}(u)^q \mathrm{~d} u
    \end{equation}
    where $\rho_{K, z}(u)=\max \{\lambda>0: z+\lambda u \in K\}$ is the radial function of $K$ with respect to $z$. When $z \in \partial K, \widetilde{V}_q(K, z)$ is defined in the way that the integral is only over those $u \in S^{n-1}$ such that $\rho_{K, z}(u)>0$. In another word,
    $$
    \widetilde{V}_q(K, z)=\frac{1}{n} \int_{\rho_{K, z}(u)>0} \rho_{K, z}(u)^q \mathrm{~d} u \text {, whenever } z \in \partial K .
    $$
    In this case, for $\mathcal{H}^{n-1}$-almost all $z \in \partial K$, we have
    $$
    \widetilde{V}_q(K, z)=\frac{1}{2 n} \int_{S^{n-1}} X_K(z, u)^q \mathrm{~d} u
    $$
    where the parallel $X$-ray of $K$ is the nonnegative function on $\mathbb{R}^n \times S^{n-1}$ defined by
    $$
    X_K(z, u)=|K \cap(z+\mathbb{R} u)|, \quad z \in \mathbb{R}^n, \quad u \in S^{n-1} .
    $$
    When $q>0$, the dual quermassintegral is the Riesz potential of the characteristic function, that is,
    $$
    \widetilde{V}_q(K, z)=\frac{q}{n} \int_K|x-z|^{q-n} \mathrm{~d} x
    $$
   Note that this immediately allows for an extension of $\widetilde{V}_q(K, \cdot)$ to $\mathbb{R}^n$. An equivalent definition via radial function can be found in \cite{LXZY2022}. By a change of variables, we obtain:
    $$
    \widetilde{V}_q(K, z)=\frac{q}{n} \int_{K-z}|y|^{q-n} \mathrm{~d} y
    $$
    since when $q>0$, the integrand $|y|^{q-n}$ being locally integrable, it can be inferred that the dual quermassintegral $\widetilde{V}_q(K, z)$ is continuous in $z$. Let $K \in \mathcal{K}^n$. The $X$-ray $X_K(x, u)$ and the radial function $\rho_{K, z}(u)$ are related as follows:
    $$
    X_K(x, u)=\rho_{K, z}(u)+\rho_{K, z}(-u), \quad \text { when } \quad K \cap(x+\mathbb{R} u)=K \cap(z+\mathbb{R} u) \neq \varnothing .
    $$
    When $z \in \partial K$, then either $\rho_{K, z}(u)=0$ or $\rho_{K, z}(-u)=0$ for almost all $u \in S^{n-1}$, and thus
    $$
    X_K(z, u)=\rho_{K, z}(u), \quad \text { or } X_K(z, u)=\rho_{K, z}(-u), \quad z \in \partial K \text {, }
    $$
    for almost all $u \in S^{n-1}$. Then, the chord integral $I_q(K)$ can be represented as follows:
    $$
    I_q(K)=\frac{1}{n \omega_n} \int_{S^{n-1}} \int_{u^{\bot}} X_K(x, u)^q \mathrm{~d} x \mathrm{~d} u, \quad q \geq 0 .
    $$
    An elementary property of the functional $I_q$ is its homogeneity. If $K \in \mathcal{K}^n$ and $q \geq 0$, then
    $$
    I_q(t K)=t^{n+q-1} I_q(K),
    $$
    for $t>0$. By compactness of $K$, it is simple to see that the chord integral $I_q(K)$ is finite whenever $q \geq 0$.
    Let $K \in \mathcal{K}^n$ and $q>0,$ the chord measure $F_q(K, \cdot)$ is a finite Borel measure on $S^{n-1}$, which can be expressed as:
    $$
    F_q(K, \eta)=\frac{2 q}{\omega_n} \int_{v^{-1}(\eta)} \widetilde{V}_{q-1}(K, z) \mathrm{d} \mathcal{H}^{n-1}(z), \quad \text { for each Borel } \eta \subset S^{n-1}.
    $$
    The mapping $v_K$ is defined on $\partial K$ with respect to the $(n-1)$-dimensional Hausdorff measure almost everywhere, due to the convexity of $K$. The chord measure $F_q(K, \cdot)$ is important as it is obtained by differentiating the chord integral $I_q$ in a certain sense, as shown in \eqref{diff}. It is evident that the chord measure $F_q(K, \cdot)$ is absolutely continuous with respect to the surface area measure $S_{n-1}(K, \cdot)$. In \cite[Theorem 4.3]{LXZY2022}, it was demonstrated that:
    $$
    I_q(K)=\frac{1}{n+q-1} \int_{s^{n=1}} h_K(v) \mathrm{d} F_q(K, v)
    $$
    When $q>0$, a useful integral formula demonstrated in \cite[Lemma 5.3]{LXZY2022} is
    \begin{footnotesize}
    	$$
    	2 n \int_{\partial K} \widetilde{V}_{q-1}(K, z) g\left(v_K(z)\right) \mathrm{d} \mathcal{H}^{n-1}(z)=\int_{s^{n-1}} \int_{\partial K} X_K(z, u)^{q-1} g\left(v_K(z)\right) \mathrm{d} \mathcal{H}^{n-1}(z) \mathrm{d} u,
    	$$
    \end{footnotesize}
    for any $g \in C\left(S^{n-1}\right)$. Therefore, for each $K \in \mathcal{K}^n$, we have
    $$
    \begin{aligned}
    \int_{S^{n-1}} g(v) \mathrm{d} F_q(K, v) & =\frac{q}{n \omega_n} \int_{S^{n-1}} \int_{\partial K} X_K(z, u)^{q-1} g\left(v_K(z)\right) \mathrm{d} \mathcal{H}^{n-1}(z) \mathrm{d} u \\
    & =\frac{q}{n \omega_n} \int_{S^{n-1}} \int_{S^{n-1}} X_K\left(\rho_K(w) w, u\right)^{q-1} h_K\left(\alpha_K(w)\right)^{-1}\\
    &\quad \rho_K(w)^n g\left(\alpha_K(w)\right) \mathrm{d} w \mathrm{~d} u .
    \end{aligned}
    $$
    Here, we denote $\rho_K=\rho_{K, o}$. For each $p \in \mathbb{R}$ and $K \in \mathcal{K}_o^n$, the $L_p$ chord measure $F_{p, q}(K, \cdot)$ is defined as follows:
    $$
    \mathrm{d} F_{p, q}(K, v)=h_K(v)^{1-p} \mathrm{~d} F_q(K, v)
    $$
    and we have an important property of $F_{p,q},$ its homogeneity, namely
    \begin{equation*}
    F_{p,q}(tK,\cdot)=t^{n+q-p-1}F_{p,q}(K,\cdot)
    \end{equation*}
    for each $t>0.$
    
    From Theorem 2.2 in \cite{XYZZ2022}, we know that if $K_i\in \mathcal{K}^{n}_o \rightarrow K_0\in \mathcal{K}^{n}_o,$ then the chord measure $F_q(K_i,\cdot)$ converges to $F_q(K,\cdot)$ weakly. Hence, one can immediately obtain that 
    $$
    F_{p,q}(K_i,\cdot) \rightarrow F_{p,q}(K,\cdot) \text{ weakly. }
    $$
    It was shown in \cite{LXZY2022} that the differential of the chord integral $I_q$ with respect to the $L_p$ Minkowski combinations leads to the $L_p$ chord measure: for $p \neq 0$,
    $$
    \left.\frac{\mathrm{d}}{\mathrm{d} t}\right|_{t=0} I_q\left(K+_p t \cdot L\right)=\frac{1}{p} \int_{S^{n-1}} h_L^p(v) \mathrm{d} F_{p, q}(K, v),
    $$
    where $K+_p t \cdot L$ is the $L_p$ Minkowski combination between $K$ and $L.$

    Since we are using the variational method to solve the $L_p$ chord Minkowski problem, the variational formula for chord integral is crucial and it is the key to tansforming the Minkowski problem into the Lagrange equation of an optimization problem. 
    \begin{thm}[Theorem 5.5 in \cite{LXZY2022}]
    	Let $q>0,$ and $\Omega$ be a compact subset of $\mathbb{S}^{n-1}$ that is not contained in any closed hemisphere. Suppose that $g:\Omega\rightarrow(0,\infty)$ is a family of continuous functions given by
    	$$
    	h_t=h_0+tg+o(t,\cdot),
    	$$
    	for each $t\in(-\delta,\delta)$ for some $\delta>0.$ Here $o(t,\cdot)\in C(\Omega)$ and $o(t,\cdot)/v$ tends to $0$ uniformly on $\Omega$ as $t\rightarrow 0.$ Let $K_t$ be the Wulff shape generated by $h_t$ and $K$ be the Wulff shape generated by $h_0.$ Then,
    	\begin{equation}\label{diff}
    	\frac{d}{dt}\big|_{t=0}I_q(K_t)=\int_{\Omega}g(v)dF_q(K,v).
    	\end{equation}
    \end{thm}
    See also in \cite[Theorem 2.1]{XYZZ2022}.
    
    Let $\mathcal{P}$ be the set of polytopes in $\mathbb{R}^n,$ and $\mathcal{P}(v_1,\cdots,v_N)$ be the subset of $\mathcal{P}$ such that $v_1,\cdots,v_N\in \mathbb{S}^{n-1}$ are in general position and every $P\in\mathcal{P}(v_1,\cdots,v_N)$ satisies
    $$
    P=\bigcap_{i=1}^{N}H^{-}(P,v_i).
    $$
    Note that, if a polytope $P\in \mathcal{P}(v_1,\cdots,v_N),$ then $P$ has at most $N$ facets, and its normals to facets are all contained in $\{v_1,\cdots,v_N\}.$ Specifically, for each $P \in \mathcal{P}\left(v_1, \ldots, v_N\right)$, the chord measure $F_q(P, \cdot)$ is entirely supported on $\left\{v_1, \ldots, v_N\right\}$. Polytopes whose facet normals are in general position have many good properties and special structure, we gather the following results .
    \begin{lem}[Lemma 4.1 \cite{Zhu14}]
    	Assuming that the unit vectors $v_1, \cdots, v_N$ are in general position and $P \in \mathcal{P}(v_1, \ldots, v_N)$, it follows that $F(P, v_i):=P\cap H(P,v_i)$ is either a point or a facet for all $1 \leq i \leq N$. Furthermore, if $n \geq 3$ and $F(P, v_i)$ is a facet, then the outer unit normals of $F(P, v_i)$ (in $H(P, v_i)$) are also in general position. 
    \end{lem}
    A key lemma obtained in Lemma 3.4 \cite{GXZ2023} reveals that if the polytopes whose facet normals are in general position get large, they will get large uniformly in all directions.
    \begin{lem}[Lemma 3.4 \cite{GXZ2023}]
    	Let $v_1,\cdots,v_N \in \mathbb{S}^{n-1}$ be in general position in dimension $n,$ and $P_i\in\mathcal{P}(v_1,\cdots,v_N).$ If outer radius $R_i$ of $P_i$ is not uniformly bound in $i,$ then its inner radius $r_i$ is not uniformly bound in $i$ either.
    \end{lem}
    We can immediately deduce that if a sequence of polytopes from $\mathcal{P}(v_1,\cdots,v_N)$ has a bounded chord integral, it also has a bounded diameter. This result is presented in the following lemma, which was proven by Xi-GZ in \cite{GXZ2023}.
    \begin{lem}\label{bound}
    	Let $v_1,\cdots,v_N \in \mathbb{S}^{n-1}$ be in general position in dimension $n,$ and $P_i\in\mathcal{P}(v_1,\cdots,v_N).$ If the $q$-th chord integral $I_q(P_i)=1,q\geqslant 0$ then the outer radius $R_i$ of $P_i$ is uniformly bounded.
    \end{lem}
    Also, in this special discrete case, the variational formula \ref{diff} turns to
    \begin{cor}
    	Let $v_1,\cdots,v_N$ be $N$ unit vectors that are in general position and $P\in\mathcal{P}(v_1,\cdots,v_N).$ Let $\delta=(\delta_1,\cdots,\delta_N)\in \mathbb{R}^{N},$ for sufficiently small $|t|,$ consider 
    	$$
    	P_t=\bigcap_{i=1}^{N}\{x\in\mathbb{R}^n:x\cdot v_i\leqslant h_P(v_i)+t\delta_i\}.
    	$$
    	Then, for $q>0,$ we have
    	$$
    	\frac{d}{dt}\big|_{t=0}I_q(P_t)=\sum_{i=1}^{N}\delta_iF_q(P,v_i).
    	$$
    \end{cor}
    
    \section{proof to theorem\ref{A}}
    Let $\mu=\sum_{i=1}^{N}\alpha_i\delta_{v_i},$ for some $\alpha_i>0$ and unit vectors $v_1,\cdots,v_N \in \mathbb{S}^{n-1}$ are in general position. Obviously, $\mu$ is a finite discrete Borel measure on $\mathbb{S}^{n-1}$ that is not concentrated in any closed hemisphere.
    
    Let $P\in\mathcal{P}(v_1,\cdots,v_N),$ define $\Phi_{p,\mu}: Int(P)\rightarrow\mathbb{R}$ by
    \begin{equation}\label{function}
    \Phi_{p,\mu}(h_P,\xi)=-\frac{1}{p}\sum_{i=1}^{N}(h_P(v_i)-\xi\cdot v_i)^p\alpha_i
    \end{equation}
    when there is no confusion what the uniderlying measure $\mu$ is, we shall write $\Phi_{p}=\Phi_{p,\mu}.$
    \begin{lem}\label{unique}
    	Let $P\in\mathcal{P}(v_1,\cdots,v_N)$ where $v_1,\cdots,v_N \in \mathbb{S}^{n-1}$ are in general position and $p<0.$ Then the minimizer of $\inf_{\xi\in P}\Phi_{p}(h_P,\xi)$ is uniquely attained at some $\xi_P\in P.$
    \end{lem}
    \begin{pf}
    	Since $\inf_{\xi\in P}\Phi_{p}(h_P,\xi)$ is invariant under all affine transformations w.r.t $h,$ we can assume $P\in \mathcal{K}^n_o.$ As $\xi\in \text{Int}P$ and $p<0,$ we have that $\Phi_{p}(h_P,\xi)\rightarrow\infty$ as $\xi\rightarrow\partial P.$ By the strict convexity of $\phi(z)=-\frac{1}{p}z^p,$ we have for $0<\lambda<1,$ and $\xi_1,\xi_2\in \text{Int}(P),$
    	\begin{footnotesize}
    		\begin{eqnarray*}
    			\Phi_{p}(h_P,\lambda\xi_1+(1-\lambda)\xi_2)&=&-\frac{1}{p}\sum_{i=1}^{N}(h_P(v_i)-(\lambda\xi_1+(1-\lambda)\xi_2)\cdot v_i)^p\alpha_i\\
    			&=&-\frac{1}{p}\sum_{i=1}^{N}(\lambda(h_P(v_i)-\xi_1\cdot v_i)+(1-\lambda)(h_P(v_i)-\xi_2\cdot v_i))^p\alpha_i\\
    			&\leqslant&-\frac{1}{p}\sum_{i=1}^{N}(\lambda(h_P(v_i)-\xi_1\cdot v_i)^p+(1-\lambda)(h_P(v_i)-\xi_2\cdot v_i)^p)\alpha_i\\
    			&=&\lambda\Phi_{p}(h_P,\xi_1)+(1-\lambda)\Phi_{p}(h_P,\xi_2),
    		\end{eqnarray*}
    	\end{footnotesize}
        with equality if and only if $\xi_1\cdot v_k=\xi_2\cdot v_k$ for all $k=1,\cdots,N,$ which implies $\xi_1=\xi_2$ since $v_1,\cdots,v_N \in \mathbb{S}^{n-1}$ are in general position. Thus $\Phi_{p}$ is strictly convex on Int$(P).$ We conclude that there exists a unique interior point $\xi_P$ such that $\Phi_{p}(h_P,\xi_P)=\min_{\xi\in\text{Int}(P)}\Phi_{p}(h_P,\xi).$  \hfill $\square$
    \end{pf}
    It is easy to see that $\xi_{\lambda P}=\lambda\xi_P,$ for $\lambda>0.$
    \begin{lem}\label{continuous}
    	If $P_i\in \mathcal{P}(v_1,\cdots,v_N)$ and $P_i$ converges to a polytope $P,$ then $\lim_{i\rightarrow\infty}\xi_{P_i}=\xi_P$ and
    	\begin{equation*}
    	\lim_{i\rightarrow\infty}\Phi_p(h_{P_i},\xi_{P_i})=\Phi_p(h_{P},\xi_{P}).
    	\end{equation*}
    \end{lem}
    \begin{pf}
    	Since $P_i$ converges to a polytope $P,$ it is simple to observe that $P\in \mathcal{P}(v_1,\cdots,v_N),$ $\xi_P\in \text{Int}P$ and $\xi_{P_i}$ are uniformly bounded in $i.$ Suppose that $\xi_{P_i}$ does not converge to $\xi_P,$ then there exists a subsequence (which we still denote as $P_i$) such that $P_{i}\rightarrow P,$ $\xi_{P_i}\rightarrow \xi_0,$ but $\xi_0\neq\xi_P.$ Note that $\xi_0\in P,$ and we have 
    	\begin{equation*}
    	\lim_{i\rightarrow\infty}\Phi_p(h_{P_i},\xi_{P_i})=\Phi_p(h_{P},\xi_{0})>\Phi_p(h_{P},\xi_{P})=\lim_{i\rightarrow\infty}\Phi_p(h_{P_i},\xi_{P}).
    	\end{equation*}
    	However, 
    	\begin{equation*}
    	\lim_{i\rightarrow\infty}\Phi_p(h_{P_i},\xi_{P})\geqslant\lim_{i\rightarrow\infty}\Phi_p(h_{P_i},\xi_{P_i})=\Phi_p(h_{P},\xi_{0}).
    	\end{equation*}
    	There comes a contradiction. Therefore, $\lim_{i\rightarrow\infty}\xi_{P_i}=\xi_P$ and thus $$\lim_{i\rightarrow\infty}\Phi_p(h_{P_i},\xi_{P_i})=\Phi_p(h_{P},\xi_{P}).$$ \hfill $\square$
    \end{pf}
    
    Consider the maximization problem
    \begin{equation}\label{optimization}
    \sup\left\{\inf_{\xi\in \text{Int}P}\Phi_{p}(h_P,\xi):P\in \mathcal{P}(v_1,\cdots,v_N) \text{ and } I_q(P)=1 \right\} .
    \end{equation}
    \begin{rmk}
    	The analogous extreme problem was studied by Chou and Wang \cite{WZ06}, Zhu \cite{Zhu14}, Xi-GZ \cite{GXZ2023} in the polytopal case.
    \end{rmk}
    \begin{lem}\label{equation}
    	If there exists a $P\in \mathcal{P}(v_1,\cdots,v_N)$ with $\xi_P=o$ and $I_q(P)=1$ such that
    	$$
    	\Phi_{p}(h_P,o)=\sup\left\{\inf_{\xi\in \text{Int}P}\Phi_{p}(h_P,\xi):P\in \mathcal{P}(v_1,\cdots,v_N) \text{ and } I_q(P)=1 \right\},
    	$$
    	then, there exists a polytope $Q\in \mathcal{P}(v_1,\cdots,v_N)$ containing the origin in its interior such that
    	$$ F_{p,q}(Q,\cdot)=\mu.$$
    \end{lem}
    \begin{pf}
    	Let $\delta=(\delta_1,\cdots,\delta_N)\in \mathbb{R}^N,$ choose $|t|$ small enough so that the polytope $Q_t$  defined by
    	$$ Q_t=\bigcap_{i=1}^{N}\left\{x:x\cdot v_i\leqslant h(P,v_i)+t\delta_i \right\}
    	$$
        contains the origin in its interior. Let $\lambda(t)=I_q(Q_t)^{-\frac{1}{n+q-1}},$ then $\lambda(t)Q_t\in \mathcal{P}(v_1,\cdots,v_N)$ with $I_q(\lambda(t)Q_t)=1$ and 
    	\begin{equation}
    	\lambda'(0)=-\frac{1}{n+q-1}\sum_{i=1}^{N}\delta_iF_q(P,v_i).
    	\end{equation}\label{lambda}
    	Let $\xi(t)=\xi_{\lambda(t)Q_t}=\lambda(t)\xi_{Q_t}$ and $$\Psi_p(t)=\Phi_{p}(h_{\lambda(t)Q_t},\xi(t)).
    	$$
    	Since $\xi(t)$ minimizes $\inf_{\xi\in \text{Int}\lambda(t)Q_t}\Phi_{p}(h_{\lambda(t)Q_t},\xi),$ we have 
    	\begin{equation}\label{xi}
    	0=\sum_{i=1}^{N}(h_{\lambda(t)Q_t}(v_i)-\xi(t)\cdot v_i)^{p-1}\alpha_iv_i.
    	\end{equation}
    	In particular, at $t=0,$ we have $0=\sum_{i=1}^{N}h^{p-1}_P(v_i)\alpha_iv_i.$ Set 
    	$$
    	F_p(t,\xi)=\sum_{i=1}^{N}(h_{\lambda(t)Q_t}(v_i)-\xi\cdot v_i)^{p-1}\alpha_iv_i.
    	$$
    	From \eqref{xi}, we know $F_p(t,\xi(t))=0.$ By a direct computation, the Jacobian with respect to $\xi$ of $F_p$ at $t=0$ and $\xi=o$ is
    	$$
    	\frac{\partial F_p}{\partial \xi}\big|_{0,o}=(1-p)\sum_{i=1}^{N}h^{p-2}_P(v_i)\alpha_i v_i\otimes v_i.
    	$$
    	From the fact that $v_1,\cdots,v_N \in \mathbb{S}^{n-1}$ are in general position, the Jacobian $\frac{\partial F_p}{\partial \xi}$ is positive-definite at $t=0$ and $\xi=o.$ By implicit function theorem, $\xi'(0)$ exists.\\
    	Since $h_{P}$ is a maximizer, we have $0=\Psi_p'(0),$ namely
    	$$
    	0=-\lambda'(0)\left(\sum_{i=1}^{N}h^{p}_P(v_i)\alpha_i\right)-\sum_{i=1}^{N}h^{p-1}_P(v_i)\alpha_i\delta_i+\xi'(0)\sum_{i=1}^{N}h^{p-1}_P(v_i)\alpha_iv_i.
    	$$
    	From \eqref{xi} and \eqref{lambda}, we have
    	$$
    	0=\frac{1}{n+q-1}\left(\sum_{i=1}^{N}h^{p}_P(v_i)\alpha_i\right)\sum_{i=1}^{N}\delta_iF_q(P,v_i)-\sum_{i=1}^{N}h^{p-1}_P(v_i)\alpha_i\delta_i.
    	$$
    	From the arbitrariness of $\delta$ we conclude that 
    	$$
    	F_{p,q}(P,\cdot)=\frac{n+q-1}{(-p)\Phi_p(h_P,o)}\mu(\cdot).
    	$$
    	The desired result immediately follows from the fact that $F_{p,q}(K,\cdot)$ is homogeneous of degree $n+q-1-p\neq0$ in $K.$ Let $Q=cP,$ where \begin{scriptsize}
    		$c=\left(\frac{(-p)\Phi_p(h_P,o)}{n+q-1}\right)^{\frac{1}{n+q-p-1}} .$
    	\end{scriptsize}  \hfill $\square$   
    \end{pf}
    The following lemma shows the existence of a solution to the maximization problem for the functional $\Phi_{p}(h_P,\xi_P)$ and by lemma \ref{equation} we know that it leads a solution to the $L_P$ chord Minkowski problem.
    \begin{thm}
    	Let $p<0,q>0,$ and $\mu=\sum_{i=1}^{N}\alpha_i\delta_{v_i},$ for some $\alpha_i>0$ and unit vectors $v_1,\cdots,v_N \in \mathbb{S}^{n-1}$ are in general position. Then there exists a polytope $P$ containing the origin in its interior such that
    	\begin{equation*}
    	F_{p,q}(K,\cdot)=\mu.
    	\end{equation*}
    \end{thm}
    \begin{pf}
    	Let $P_i\in\mathcal{P}(v_1,\cdots,v_N)$ be a maximizing sequence; that is $I_q(P_i)=1$ and 
    	$$
    	\lim_{i\rightarrow\infty}\Phi_p(h_{P_i},\xi_{P_i})=\sup\left\{\Phi_{p}(h_P,\xi_P):P\in \mathcal{P}(v_1,\cdots,v_N) \text{ and } I_q(P)=1 \right\} .
        $$
        By translation invariance of $I_q$ and $\Phi_{p},$ for $P,Q\in\mathcal{P}(v_1,\cdots,v_N),$ if there exists an $x\in\mathbb{R}^n$ such that $P=Q+x,$ then $$ \Phi_p(h_P,\xi_P)=\Phi_p(h_Q,\xi_Q).$$
        Thus, we can assume $\xi_{P_i}=o.$ From lemma\ref{bound} we have $P_i$ is uniformly bounded. Hence, we can choose a subsequence (which we still denote as $P_i$) such that $P_i\rightarrow P$ for some $P\in\mathcal{P}(v_1,\cdots,v_N)$ by the Blaschke selection theorem. By continuity of $I_q$, we have $I_q(P)=1$ which implies $P$ is non-degenerated. From lemma\ref{continuous} we have $\xi_P=o\in \text{Int}P$ and 
        $$
        \Phi_{p}(h_P,\xi_P)=\lim_{i\rightarrow\infty}\Phi_p(h_{P_i},\xi_{P_i})=\sup_{P\in \mathcal{P}(v_1,\cdots,v_N)}\left\{\Phi_{p}(h_P,\xi_P):  I_q(P)=1 \right\} .
        $$
        The desired result now immediately follows from lemma\ref{equation}. \hfill $\square$
    \end{pf}
    \section{proof to the theorem\ref{B}}
    Let $\mu$ be a finite Borel measure on $\mathbb{S}^{n-1}$ that is not concentrated in any closed hemisphere with a density fanction $\frac{1}{\Lambda}<f<\Lambda.$\\
    We consider a similar functional as in the discrete case, but, we have to consider instead an approximation problem first. This is because we need the positivity of the maximizer of a maximization problem \eqref{max}, however, when $p \in (1-n,0),$ the maximizer may fail to be positive. For $\varepsilon \in (0,1/4)$ small, let $\phi_{\epsilon}$ be a positive, strictly convex, monotone decreasing function on $(0, \infty)$ such that
    $$
    \phi_{\epsilon}(z) \begin{cases}=-\frac{1}{p} z^{p}, & z \geqslant 1, \\ \geqslant \frac{1}{n-1} z^{-n+1}, & 0<z<\varepsilon \\ \leqslant -\frac{1}{p'} z^{p^{\prime}}, & 0<z<\frac{1}{4}\end{cases}
    $$
    where $p'\in(-n,-n+1]$ is a fixed constant, For $p\in(-n,-n+1],$ we can take
    $$
    \phi_{\epsilon}(z)=-\frac{1}{p} z^{p}
    $$
    for all $z>0.$ Let $\Phi_p(h)=\int_{\mathbb{S}^{n-1}}f\phi_{\epsilon}(h).$
    From a similar proof as in lemma\ref{unique} and lemma\ref{continuous}, we can conclude that there exists a unique $\xi_h\in \text{Int }K_h$ which attains the infimum of  $\inf\left\{\Phi_p(h-\xi\cdot x):\xi\in K\right\},$ and we also have the continuity of $\xi_h$ in $h.$ In cases where there is no confusion, we may write $\xi_{h}=\xi_{K}$. We conclude the results in the following lemma.
    \begin{lem}
    	Let $f,\phi_{\epsilon}, \Phi_{p}$ be defined as above, $K\in \mathcal{K}_{o}^{n},$ and $-n<p<0.$ Then the minimizer of $\inf\left\{\Phi_p(h-\xi\cdot x):\xi\in K\right\}$ is uniquely attained at some $\xi_{K}\in K.$ Moreover, if $\{K_i\}\subset \mathcal{K}_{o}^{n}$ and $K_i$ converges to a convex body $K,$ then $\lim_{i\rightarrow\infty}\xi_{K_i}=\xi_{K}$ and 
    	$$
    	\lim_{i\rightarrow\infty}\Phi_{p}(h_{K_i}-\xi_{K_i}\cdot x)=\Phi_{p}(h_{K}-\xi_{K}\cdot x).
    	$$
    \end{lem}
    \begin{pf}
    	The existence and uniqueness of $\xi_K$ can be deduced directly from the strict convexity of $\phi_{\epsilon}$ and the fact that $\xi \in \text{Int }K,p<0.$ 
    	
    	When $K_i \in \mathcal{K}_{o}^{n}$ converges to a convex body $K,$ $K_i$ is bounded and so is $\xi_{K_i}.$ Then there exists a subsequence $K_{i_j}$ of $K_i$ satisfying $K_{i_j}\rightarrow K,$ $\xi_{K_{i_j}}\rightarrow \xi$ for some $\xi\in K.$ Since 
    	$$
    	\lim_{j\rightarrow\infty}\Phi_{p}(h_{K_{i_j}}-\xi_{K_{i_j}}\cdot x)=\lim_{j\rightarrow\infty}\Phi_{p}(h_{K}-\xi\cdot x)=\infty
    	$$
    	if $\xi$ is a boundary point of $K.$ But 
    	$$
    	\Phi_{p}(h_{K_{i_j}}-\xi_{K_{i_j}}\cdot x)\leqslant \Phi_{p}(h_{K_{i_j}}-\xi_{K}\cdot x)< \int_{\mathbb{S}^{n-1}}f\phi_{\epsilon}(\frac{a}{2})
    	$$
    	whenever $j$ big enough and $a=\min_{v\in \mathbb{S}^{n-1}}\{h_{K}-\xi_K\cdot x\}>0.$
    	Hence, we know that $\xi$ is an interior point of $K.$ We only need to show that $\xi=\xi_K.$ If not, we have
    	\begin{eqnarray*}
    	\lim_{j\rightarrow\infty}\Phi_{p}(h_{K_{i_j}}-\xi_{K_{i_j}}\cdot x)&=&\Phi_{p}(h_{K}-\xi\cdot x)\\
    	&>&\Phi_{p}(h_{K}-\xi_K \cdot x)\\
    	&=&\lim_{j\rightarrow\infty}\Phi_{p}(h_{K_{i_j}}-\xi_{K }\cdot x).
    	\end{eqnarray*}
        This will cause a contradiction with $\Phi_{p}(h_{K_{i_j}}-\xi_{K_{i_j}}\cdot x)\leqslant \Phi_{p}(h_{K_{i_j}}-\xi_{K}\cdot x).$ Hence we have $\lim_{i\rightarrow\infty}\xi_{K_i}=\xi_{K}$ and 
        $$
        \lim_{i\rightarrow\infty}\Phi_{p}(h_{K_i}-\xi_{K_i}\cdot x)=\Phi_{p}(h_{K}-\xi_{K}\cdot x).
        $$ \hfill $\square$
    \end{pf}
    Next, we consider the maximization problem
    \begin{equation}\label{max}
    \sup_{h\in\mathcal{S}^{+}}\left\{\inf_{\xi\in K_h}\Phi_p(h-\xi\cdot x): I_q(K_h)=1\right\}.
    \end{equation}
    Before we begin the proof of the main theorem, it is important to recall two key inequalities that will be utilized in the subsequent analysis.
    One is the Blaschke-Santalo inequality
    $$
    \sup_{h\in\mathcal{S}^{+}}\inf_{\xi\in K}V(h)\int_{\mathbb{S}^{n-1}}\frac{1}{(h-\xi\cdot x)^{n}}\leqslant \frac{\omega_{n-1}^2}{n}
    $$
    where $K=K_h$ denotes the convex body determined by $h$, and the infimum is taken over all $\xi$ satisfying $h-\xi\cdot x\in \mathcal{S}^{+}$. It is worth noting that the left-hand side of this inequality is invariant under all affine transformations. The equality holds if and only if $K$ is a centered ellipsoid. Another significant inequality is presented in the following lemma, which was obtained in \cite{LXZY2022}:
    \begin{lem}
    	If $K\in \mathcal{K}^n_o$ and $1\leqslant r<s,$ then 
    	$$
    	I_r(K)\leqslant c(s,r)V(K)^{1-\frac{r-1}{s-1}}I_s(K)^{\frac{r-1}{s-1}},
    	$$ 
    	with $c(s,r)=rs^{-\frac{r-1}{s-1}}.$
    \end{lem}
    This inequality follows from a simple argument using jensen's inequality. When we choose $1\leqslant r=q < n+1,$ we have 
    \begin{eqnarray}\label{volume}
    I_q(K)&\leqslant& c(n+1,q)V(K)^{1-\frac{q-1}{n+1-1}}I_{n+1}(K)^{\frac{q-1}{n+1-1}}\\
    &=&c(n+1,q)V(K)^{\frac{n-q+1}{n}}\{\frac{n+1}{\omega_n}V(K)^2\}^{\frac{q-1}{n}}\\
    &=& q\omega_n^{-\frac{q-1}{n}}V(K)^{\frac{n+q-1}{n}},
    \end{eqnarray} 
    which is important for the following analysis.
    %In particular, we have 
   % \begin{equation}\label{volume}
    %I_q(K)\leqslant I_q(B_n)(\frac{V(K)}{\omega_n})^{\frac{n+q-1}{n}},q\geqslant 1.
   % \end{equation}
   % as a special case of the lemma.
    %    Denote $c=\sup_{h\in\mathcal{S}^{+}}\left\{\inf_{\xi\in K_h}\Phi_p(h-\xi\cdot x): I_q(K_h)=1\right\},$ we claim that there exi
    
    Frist, we demonstrate that the maximization problem does have a solution. The key point is to prove the convergence of a maximizing sequence of convex bodies is a convex body. That is to prove compactness and non-degeneracy. The assumption $I_q(K)=1$ can make sure the non-degeneracy, so we only need to prove the compactness. Here we let $n+1>q\geqslant 1$ and we can apply \eqref{volume} to obtain a uniform lower bound of volume, then we can follow the same strategy as in \cite[lemma 5.1]{WZ06}. 
    \begin{lem}\label{exi2}
    	The maximization problem \eqref{max} has a solution $h_{\epsilon}\in \mathcal{S}^{+}$ with $\xi_{h_{\epsilon}}=o$.
    \end{lem}
    \begin{pf}
    	Let $h_i\in \mathcal{S}^+$ be a maximizing sequence; that is $I_q(K_i)=1$ and 
    	$$
    	\lim_{i\rightarrow\infty}\Phi_p(h_i-\xi_{h_i})=\sup_{h\in\mathcal{S}^{+}}\left\{\inf_{\xi\in K_h}\Phi_p(h-\xi\cdot x): I_q(K_h)=1\right\}
    	$$
    	where $K_i$ is the convex body uniquely determined by $h_i.$ We claim that the diameter of $K_i$ is uniformly bounded by some constant $C>0$ independent of $\epsilon.$ Denote $d_i$ as the diameter of $K_i,$ we need to show
    	\begin{equation}\label{diam2}
    	d_i\leqslant C.
    	\end{equation}
    	Suppose on the contrary that $d_i\rightarrow\infty$ as $i\rightarrow\infty.$ Since $\inf_{\xi\in K_h}\Phi_p(h-\xi\cdot x)$ and $I_q(K_h)$ are invariant under any translation of $K_h,$ we may assume $K_i$ is centered at the origin and so as the minimium ellipsoid $E_i$ of $K_i$ by John' lemma. Then we have $c(n)E_i\subset K_i\subset E_i.$ Let $\delta\in(0,1/4)$ be a fixed small constant, $
    	S_1=\mathbb{S}^{n-1}\cap\{h_{E_i}<\delta\},$ $S_2=\mathbb{S}^{n-1}\cap\{\delta \leqslant h_{E_i} \leqslant 1/\delta\},$ $S_3=\mathbb{S}^{n-1}\cap\{h_{E_i}>1/\delta\}.$ Obviously, we have $\mathbb{S}^{n-1}=S_1\cup S_2\cup S_3.$ Since $c(n)h_{E_i}\leqslant h_i \leqslant h_{E_i},$ then by the monotoniciti of $\phi_{\epsilon}$
    	\begin{eqnarray*}
    		\inf_{\xi\in K_i}\Phi_p(h_i-\xi\cdot x)&\leqslant& \Phi_p(h_i)\\
    		&\leqslant& \Phi_p(c(n)h_{E_i})\\
    		&=& \int_{\mathbb{S}^{n-1}}f\phi_{\epsilon}(c(n)h_{E_i})\\
    		&=& \sum_{j=1}^{3}\int_{S_j}f\phi_{\epsilon}(c(n)h_{E_i}).
    	\end{eqnarray*}
    	By definition of the function $\phi_{\epsilon}$ and the set $S_1$ we have
    	\begin{eqnarray*}
    		\int_{S_1}f\phi_{\epsilon}(c(n)h_{E_i})&\leqslant& \int_{S_1}\frac{C}{h_{E_i}^{-p'}}\\
    		&\leqslant& C\left(\int\frac{1}{h_{E_i}^{n}}\right)^{\frac{-p'}{n}}|S_1|^{\frac{n+p'}{n}}\\
    		&\leqslant& C|S_1|^{\frac{n+p'}{n}}\\
    		&\rightarrow& 0,
    	\end{eqnarray*}
    	for any fixed $\delta,$ the second is from the H$\ddot{o}$lder inequality and the third one is by the Blaschke-Santalo inequality and \eqref{volume}. Indeed, since $I_q(K_i)=1,$ by \eqref{volume} we have $V(E_i)\geqslant V(K_i)\geqslant c$ for some constant $c$ independent on $i.$ The last one is from the assuption that $d_i\rightarrow\infty,$ then $d_{E_i}\rightarrow\infty,$ but we have $1\leqslant I_{q}(E_i)\leqslant c(n)^{-(n+q-1)},$ which implies that the inner radius can be as small as one hopes as $i\rightarrow\infty,$ and thus $|S_1|,|S_2|\rightarrow 0$ as $i\rightarrow\infty.$ Hence, we also have $\int_{S_2}f\phi_{\epsilon}(c(n)h_{E_i})\rightarrow 0,$ then 
    	\begin{eqnarray}\label{critical}
    	\int_{\mathbb{S}_{n-1}}f\phi_{\epsilon}(c(n)h_{E_i})&=& o(1)+\int_{S_3}f\phi_{\epsilon}(c(n)h_{E_i})\\
    	&\leqslant& o(1)+C\phi_{\epsilon}\left(\frac{c(n)}{\delta}\right) 
    	\end{eqnarray}
    	where $C$ is independent of $i,\delta,\epsilon$ for each $i\rightarrow \infty.$ Denote $$c=\sup_{h\in\mathcal{S}^{+}}\left\{\inf_{\xi\in K_h}\Phi_p(h-\xi\cdot x): I_q(K_h)=1\right\},$$ \eqref{critical} implies 
    	\begin{equation}\label{up}
    	c\leqslant C\phi_{\epsilon}(\delta^{-1}).
    	\end{equation}
    	On the other hand, taking $h\equiv r$ for some $r<1$ such that $I_q(B_r)=1.$ Indeed, since $I_q(B_r)=r^{n+q-1}I_q(B_1),$ we can choose $r=\frac{\omega_q}{2^q\omega_n\omega_{n+q-1}}.$ Hence, we have
    	\begin{eqnarray}
    	c&\geqslant& \inf_{\xi\in B_r}\Phi_p(r-\xi\cdot x)\\
    	&\geqslant& \int_{\mathbb{S}^{n-1}}f\phi_{\epsilon}(2)\\
    	&\geqslant& C_1\phi(2)>0,
    	\end{eqnarray}
    	by the monotone decreasing condition of $\phi$ and $|\xi\cdot x|\leqslant1,$ and we can use $\phi_{\epsilon}(z)\leqslant C+z^{-n}$ for $z>0$ and the Blaschke-Santalo inequality to abtain the boundness from above, that is
    	$$
    	\inf_{\xi\in K_h}\Phi_p(h-\xi\cdot x)\leqslant C+\inf_{\xi\in K_h}\int\frac{f}{(h-\xi\cdot x)^{n}}\leqslant C
    	$$
    	for all $h\in\mathcal{S}^{+},I_q(K_h)=1.$ Hence, we have $c_1<c\leqslant C$ for some constant $c_1,C$ independent of $\epsilon,$ this inequality is in conflict with \eqref{up} since $\phi_{\epsilon}(\delta^{-1})$ tends to $0.$
    	
    	Therefore, $\{h_i\}$ converges to a maximizer $h_{\epsilon}$ of \eqref{max} up to a subsequence, and the diameter of $K_{h_{\epsilon}}$ is also bounded by a constant $C>0$ independent of $\epsilon.$ By a translation we may assume that $\xi_{h_{\epsilon}}=o,$ which implies that $h_{\epsilon}>0$ on $\mathbb{S}^{n-1}.$ \hfill $\square$
    \end{pf}
    \begin{lem}
    	The maximizer $h_{\epsilon}$ obtained in Lemma\ref{exi2} satisfies
    	\begin{equation}\label{c0 estimate}
    	c\leqslant h_{\epsilon}\leqslant C
    	\end{equation}
    	for some constants $c,C$ that are independent of $\epsilon$.
    \end{lem}
    \begin{pf}
    	Let $K_{\epsilon}$ be the convex body corresponding to $h_{\epsilon},$ that is to say, $h_{\epsilon}$ is the support function of $K_{\epsilon}.$ We have already established a uniform lower bound on the volume and the right-hand inequality $h_{\epsilon}\leqslant C$ through the proof of Lemma \ref{exi2}. Therefore, we only need to demonstrate that $h_{\epsilon}$ is bounded from below. Suppose that $h_{\epsilon}(v_{\epsilon})\rightarrow 0$ for some $v_{\epsilon}\in \mathbb{S}^{n-1},$ as $\epsilon \rightarrow 0$.Due to the compactness of $\mathbb{S}^{n-1}$, there exists $v\in \mathbb{S}^{n-1}$ such that $h_{\epsilon}(v)\rightarrow 0$ as $\epsilon \rightarrow 0$ (up to a subsequence of $h_{\epsilon}$).
    	
        Since there exist $R$ and $\hat{c}$, both independent of $\epsilon$, such that 
        $$
        K_{\epsilon}\subset B_R(O) \text{ and } V(K_{\epsilon})\geqslant \hat{c},
        $$
        there exists $r>0$ that depends only on $n$, $R$, and $c>0$, such that 
        $$
        B_r(p_{\epsilon})\subset K_{\epsilon},
        $$ 
        where $p_{\epsilon}$ is the barycenter of $K_{\epsilon}$. Note that $h_{\epsilon}$ is a maximizer of the maximization problem \eqref{max} with $\xi_{\epsilon}=o.$ Hence, by assumption we have that $\xi_{\epsilon} \rightarrow \partial K_{\epsilon}$ as $\epsilon \rightarrow 0$. Hence, we have 
        $$
        \Phi_{p}(h_{\epsilon}-\xi_{\epsilon}\cdot x)=\Phi_{p}(h_{\epsilon})\rightarrow \infty,
        $$
        as $\epsilon \rightarrow 0$. This contradicts the fact that 
        $$
        C\geqslant c=\sup_{h\in\mathcal{S}^{+}}\left\{\inf_{\xi\in K_h}\Phi_p(h-\xi\cdot x): I_q(K_h)=1\right\}.
        $$ Therefore, we have the uniformly $C^0$ estimate $c\leqslant h_{\epsilon}\leqslant C$. \hfill $\square$
    \end{pf}
    With the uniform boundness of $h_{\epsilon}$, we then can deduce the bounds of the $q-1$ th dual quermassintegral $\widetilde{V}_{q-1}([h_{\epsilon}],\bar{\nabla} h_{\epsilon})$ from both sides independent of $\epsilon$ and the point $\bar{\nabla} h_{\epsilon}.$ 
    \begin{lem}
    	 Let $h_{\epsilon}$ be the maximizer obtained in Lemma\ref{exi2}, $K_{\epsilon}$ be the convex body corresponding to $h_{\epsilon}$ and $q\geqslant 1,z\in \partial K_{\epsilon}$, then there exists two constants $c,C$ such that
    	\begin{equation}\label{quermass}
    	c\leqslant \widetilde{V}_{q-1}(K_{\epsilon},z)\leqslant C
    	\end{equation}
    	where $c,C$ are independent with $\epsilon$ and $z.$
    \end{lem}
    \begin{pf}
    	By \eqref{c0 estimate}, we have that 
    	$$
    	B_c(o)\subset K_{\epsilon}\subset B_C(O)
    	$$
    	where $c,C$ are the same constants with that in \eqref{c0 estimate}.
    	Since 
    	$$
    	\widetilde{V}_{q-1}(K_{\epsilon}, z)=\frac{1}{2 n} \int_{S^{n-1}} X_{K_{\epsilon}}(z, u)^{q-1} d u
    	$$
    	where the parallel $X$-ray of $K_{\epsilon}$ is the nonnegative function on $\mathbb{R}^n \times S^{n-1}$ defined by
    	$$
    	X_{K_{\epsilon}}(z, u)=|K_{\epsilon} \cap(z+\mathbb{R} u)|, \quad z \in \mathbb{R}^n, \quad u \in S^{n-1} .
    	$$
    	We can choose a Borel set $Z\subset \mathbb{S}^{n-1}$ satisfying $|Z|\geqslant c,$ here, $c$ is a universal constant independent with $z$ and $\epsilon,$ and $\forall u\in Z$ we have $$X_{K_{\epsilon}}(z, u)\geqslant \tilde{c} $$ for some uniform constant $\tilde{c}$. Indeed, in dimension 2, since $K_{\epsilon}$ is pinched between two bounded balls, $\forall z\in \partial K_{\epsilon},$ there is a Borel set $Z\subset \mathbb{S}^{1}$ with $\arcsin(\frac{\sqrt{3}c}{2C})\leqslant |Z|\leqslant \frac{2\pi}{3}$ such that $u\in Z,X_{K_{\epsilon}}(z, u)\geqslant c,$ here, the constant $c(C)$ is the radius of the inner(outer) ball of $K_{\epsilon}.$ The higher dimensional case is analogous, just do a rotation. Hence 
    	$$
    	\begin{aligned}
    	\widetilde{V}_q(K_{\epsilon}, z) &=\frac{1}{2 n} \int_{S^{n-1}} X_{K_{\epsilon}}(z, u)^q d u \\
    	&\geqslant \frac{1}{2 n} \int_{Z} X_{K_{\epsilon}}(z, u)^q d u\\
    	&\geqslant \frac{1}{2 n} \int_{Z} c^q d u\\
    	&\geqslant c_1 |Z|\\
    	&\geqslant c_2.
    	\end{aligned}
    	$$ 
    	On the other hand, since 
    	$$
    	X_{K_{\epsilon}}(z, u)\leqslant 2\max h_{\epsilon}\leqslant 2C,
    	$$
    	which implies
    	$$
    	\begin{aligned}
    	\widetilde{V}_q(K_{\epsilon}, z) &=\frac{1}{2 n} \int_{S^{n-1}} X_{K_{\epsilon}}(z, u)^q d u \\
    	&\leqslant \frac{1}{2 n} \int_{S^{n-1}} (2C)^q d u\\
    	&\leqslant C_3.
    	\end{aligned}
    	$$
        \hfill $\square$
    \end{pf}
    We are in a position to construct a variation of the optimization problem and obtain its Lagrange equation. We show that maximizer of the optimization problem is a generalied solution of a Monge-Amp\'{e}re type equation:
    \begin{thm}\label{approeq}
    	Let $h_{\epsilon}$ be the maximizer in lemma \ref*{exi2}. Then $h_{\epsilon}$ is a generalized solution of
    	$$
    	\mbox{det}(\nabla^2h_{\epsilon}+h_{\epsilon}I) =-\frac{f\phi_{\epsilon}'}{C_{\epsilon} \widetilde{V}_{q-1}([h_{\epsilon}],\bar{\nabla} h_{\epsilon})},
    	$$
    	where $C_{\epsilon}=-\frac{1}{n+q-1}\int f\phi_{\epsilon}'h_{\epsilon}.$
    \end{thm}
    \begin{pf}
    	Since $h_{\epsilon}$ is the maximizer obtained from lemma \ref{exi2}, $h_{\epsilon}>0.$ Suppose $g\in C(\mathbb{S}^{n-1}),$ let
    	$$
    	K_t=\left\{x\in\mathbb{R}^n:x\cdot u\leqslant(h_{\epsilon}+tg)(u)\right\} 
    	$$
    	for sufficiently small $\delta>0$ such that $|t|<\delta$ then $h_{\epsilon}+tg>0.$ Let $h_t$ be the support function of $K_t,$ Note that $h_0=h_{\epsilon},K_0=K.$ 
    	Let $\lambda(t)=I_q(K_t)^{-\frac{1}{n+q-1}},$ then $I_q(\lambda(t)K_t)=1,$ and $\lambda'(0)=-\frac{1}{n+q-1}\int_{\mathbb{S}^{n-1}}g(v)dF_q(K,v).$ Denote 
    	$$
    	\Psi_p(t)=\Phi_{p}(\lambda(t)h_t-\lambda(t)\xi_t \cdot x),
    	$$
    	where $\xi_t$ is the minimizer of $\inf_{\xi\in K_t}\Phi_{p}(h_t-\xi \cdot x).$ Based on the assumption about $\phi_{\epsilon}$, it can be concluded that $\xi$ is Lipschitz continuous. Without loss of generality, let us assume
    	$$
    	\lim _{t_j \rightarrow 0} \frac{\xi\left(t_j\right)-\xi(0)}{t_j}=\xi'(0) ,
    	$$
    	this assumption is valid, and we can follow the same procedure as in the proof of Lemma \ref{equation} to show that $\xi'(0)$ is well-defined. As $h_{\epsilon}$ is a maximizer, the function $t\mapsto \Psi_p(t)$ attains its maximum at $t=0$. However, $ h_t$ may not be differentiable at $t=0.$ Let
    	$$
    	\psi_p(t)=\Phi_{p}(\lambda(t)(h_{\epsilon}+tg)-\lambda(t)\xi_t \cdot x).
    	$$
    	Since $K_t$ is the Wulff shape of $h_{\epsilon}+tg$, we have $h_t\leqslant h_{\epsilon}+tg.$ By the monotone decreasing condition of $\phi$ and $f\geqslant \frac{1}{\Lambda}, \lambda(t)>0,$ we have 
    	$$
    	\lambda(t)h_t-\lambda(t)\xi_t \cdot x\leqslant \lambda(t)(h+tg)-\lambda(t)\xi_t \cdot x
    	$$
    	and thus
    	$$
    	\Psi_p(0)\geqslant \Psi_p(t)\geqslant \psi_p(t).
    	$$
    	Since $\Psi_p(0)=\psi_p(0),$ $\psi_p(t)$ also attains its maximum at $t=0.$ Hence we have
    	$$
    	\int f \phi^{\prime}_{\epsilon}\left(\lambda^{\prime}(0) h_{\epsilon}+g-\xi'(0)\cdot x\right) = 0 .
    	$$
    	Since the infimum of $\int_{S^n} f \phi_{\epsilon}(h_{\epsilon}-\xi \cdot x)$ is attained at $\xi=0$. We have
    	$$
    	\int_{S^{n-1}} f \phi^{\prime}_{\epsilon} x_i=0, \quad i=1, \ldots, n .
    	$$
    	Therefore,
    	$$
    	\int_{S^{n-1}} f \phi^{\prime}_{\epsilon}\left(\lambda^{\prime}(0) h_{\epsilon}+g\right) = 0 .
    	$$
    	It follows that
    	$$
    	\int f \phi^{\prime}_{\epsilon} g = \frac{1}{n+q-1}\int f\phi^{\prime}_{\epsilon}h_{\epsilon}\int g d F_q(K,\cdot) .
    	$$
    	Since $g \in C\left(S^n\right)$ is arbitrary, we conclude that 
    	$$
    	\mbox{det}(\nabla^2h_{\epsilon}+h_{\epsilon}I) =-\frac{f\phi^{\prime}_{\epsilon}}{C_{\epsilon} \widetilde{V}_{q-1}([h_{\epsilon}],\bar{\nabla} h_{\epsilon})},
    	$$
    	where $C_{\epsilon}=-\frac{2q}{(n+q-1)\omega_n}\int f\phi^{\prime}_{\epsilon}h_{\epsilon}.$ \hfill $\square$
    \end{pf}
    Now, we can prove Theorem \ref{B}. Let $\phi_{\epsilon}$ be as above. We may further assume that $\phi_{\epsilon}$ satisfies
    $$
    \lim _{\epsilon \rightarrow 0} \phi_{\epsilon}(z)=\phi_0(z)=-\frac{1}{p} z^{p}, 
    $$
    and
    $$
    \lim _{\epsilon \rightarrow 0} \phi_{\epsilon}^{\prime}(z)=-z^{p-1},
    $$
    uniformly on every compact subset of $(0, \infty)$. Let $h_{\epsilon}$ be the maximizer in Lemma \ref{exi2}. By Theorem \ref{approeq}, $h_{\epsilon}$ satisfies
    \begin{equation}\label{epsilonMA}
    \mbox{det}(\nabla^2h_{\epsilon}+h_{\epsilon}I) =-\frac{f\phi^{\prime}_{\epsilon}}{C_{\epsilon}\widetilde{V}_{q-1}([h_{\epsilon}],\bar{\nabla} h_{\epsilon})},
    \end{equation}
    where $C_{\epsilon}=-\frac{2q}{(n+q-1)\omega_n}\int f\phi^{\prime}_{\epsilon}h_{\epsilon}>0.$
%    The definition of $\phi_{\epsilon}$ implies that $h_{\epsilon}$ is positive.
    % and hence the full regularity of $h_{\epsilon}$ follows from the general theory.
    
    Since $c \leqslant h_{\epsilon}\leqslant C$ for some constant independent on $\epsilon,$ by passing to subsequences, we may suppose that $h_{\epsilon} \rightarrow h_0$ and hence $C_{\epsilon} \rightarrow C_0$ as $\epsilon \rightarrow 0$. By the continuity of the $q$th chord integral we have $I_q(K_0)=1$ and $K_0$ has bounded diameter. Additionally, we know that the infimum of $\int_{S^n} f \phi_0\left(h_0-\xi \cdot x\right)$ is attained at $\xi=0$. We claim that $C_0<\infty$. If not, $C_0=\infty$, then
    $$
    \mbox{det}(\nabla^2h+hI)=0 \quad \text { on }\left\{h_0>0\right\},
    $$
    which implies that the area measure vanishes on $\left\{h_0>0\right\}$. This is impossible. On the other hand, we say $C_0>0$. Otherwise by \eqref{quermass} we will have
    $$
    \begin{aligned}
    -\int_{S^{n-1}} f \phi_{\epsilon}^{\prime}\left(h_{\epsilon}\right) & =C_{\epsilon} \int_{S^{n-1}} \mbox{det}(\nabla^2h_{\epsilon}+h_{\epsilon}I)\widetilde{V}_{q-1}([h_{\epsilon}],\bar{\nabla} h_{\epsilon}) \\
    & \leqslant C C_{\epsilon}S(K_{\epsilon}) \\
    & \longrightarrow 0 .
    \end{aligned}
    $$
    However, this is impossible too.
    Finally, we have to check that $h_0$ solves \eqref{MAeq}. From the weak continuity property for $L_p$ chord measures, we have
    $$
    F_{p,q}(K_{\epsilon},\cdot) \rightarrow F_{p,q}(K_0,\cdot) \text{ weakly. }
    $$
    This implies $h_0$ is a generalized solution of
    $$
    \mbox{det}(\nabla^2h_{0}+h_{0}I) =\frac{f(h_{0})^{p-1}}{C_{0}\widetilde{V}_{q-1}([h_{0}],\bar{\nabla} h_{0})}.
    $$
    Moreover, using the uniform bounds \eqref{c0 estimate} and \eqref{quermass}, by the regularity results in \cite{Caf 90C1}, we know the right hand side of the above equation is bounded away from zero and infinity, we know $h_{0}\in C^1.$
    After a suitable scaling $\alpha h_0$ solves \eqref{MAeq} with $\alpha =C_0^{\frac{1}{n+q-p-1}}$. The proof of Theorem \ref{B} is completed.
    
     Declaration

    Conflict of interest The authors declare that there are no conflict of interest regarding the publication of
    this paper.

\bibliographystyle{amsplain}

\begin{thebibliography}{12}
	\bibitem{Andrews (1)}
	B. Andrews, Gauss curvature flow: the fate of the rolling stones, Invent. Math. 138(1999), 151-161.
	
	\bibitem{BLYZ13}
	K.J. B$\ddot{o}$r$\ddot{o}$czky, E. Lutwak, D. Yang and G. Zhang, The logarithmic Minkowski problem, J. Amer. Math. Soc. 26(2013), 831-852.
	
	\bibitem{BHZ16}
	K.J. B$\ddot{o}$r$\ddot{o}$czky, P. Hegedűs and G. Zhu, On the discrete logarithmic Minkowski problem, Int. Math. Res. Not. IMRN 6(2016), 1807-1838.
	
	\bibitem{BLYZ12(5)}
	K.J. B$\ddot{o}$r$\ddot{o}$czky, E. Lutwak, D. Yang and G. Zhang, The log-Brunn-Minkowski inequality, Adv. Math. 231(2012), 1974-1997.
	
	\bibitem{BCD(7)}
	S. Brendle, K. Choi and P. Daskalopoulos, Asymptotic behavior of flows by powers of the Gaussian curvature, Acta Math. 219(2017), 1-16.
	
	\bibitem{Caf 90C1}
	L.A. Caffarelli, A localization property of viscosity solutions to the Monge–Ampire equation and their strict convexity, Ann. Math. 131 (1990) 129–134.
	
	\bibitem{WZ06}
	K.-S. Chou X.-J. Wang. The Lp-Minkowski problem and the Minkowski problem in centroaffine geometry. Adv. Math. 205 (2006), no. 1, 33–83. MR2254308.
	
	\bibitem{Cheng-Yau}
	S.Y. Cheng, S.T. Yau, On the regularity of the n-dimensional Minkowski problem, Comm. Pure
	Appl. Math. 20 (1977) 41–68.
	
	\bibitem{CHLL(10)}
	S. Chen, Y. Huang, Q.-R. Li and J. Liu, The Lp-Brunn-Minkowski inequality for $p < 1$, Adv. Math. 368(2020), 107166.
	
	\bibitem{CFL 2022}
	S. Chen, Y. Feng, W. Liu, Uniqueness of solutions to the logarithmic Minkowski problem in $\mathbb{R}^3$. Adv. Math. 411 (2022), part A, Paper No. 108782, 18 pp. 52A40, MR4512403.
	
	\bibitem{Firey(17)}
	W.J. Firey, Shapes of worn stones, Mathematika 21(1974), 1-11.
	
	\bibitem{Gage18}
	M. Gage, Evolving plane curves by curvature in relative geometries, Duke Math. J. 72(1993), 441-466.
	
	\bibitem{LYZH2005}
	D. Hug, E. Lutwak, D. Yang, and G. Zhang. On the Lp Minkowski problem for polytopes. Discrete Comput. Geom., 33:699–715, 2005.
	
	\bibitem{KL polytope}
	G. Károlyi, L. Lovász, Decomposition of convex polytopes into simplices, preprint.
	
	\bibitem{LXZY2022}
	E. Lutwak, D. Xi, D. Yang, and G. Zhang. Chord measure in integral geometry and their Minkowski problems. Comm. Pure Appl. Math., in press.
	
	\bibitem{L1}
	E. Lutwak, The Brunn–Minkowski–Firey theory I: mixed volumes and the Minkowski problem,
	J. Differential Geom. 38 (1993) 131–150 II: affine and geominimal surface areas, Adv. Math.
	118 (1996) 244–294.
	
	\bibitem{Pogorelov}
    A.V. Pogorelov, The Multidimensional Minkowski Problem, Wiley, New York, 1978.
    
    \bibitem{D. Ren}
    D. Ren, Topics in Integral Geometry, World Scientific, Singapore, 1994.
	
	\bibitem{Stancu0240}
	A. Stancu, The discrete planar $L_0$ Minkowski problem, Adv. Math. 167(2002), 160-174.
	
	\bibitem{Stancu41}
	A. Stancu, On the number of solutions to the discrete two-dimensional $L_0$ Minkowski problem, Adv. Math. 180(2003), 290-323.
	
	\bibitem{Schneider}
	R. Schneider. Convex bodies: the Brunn-Minkowski theory, volume 151 of Encyclopedia of Mathematics and
	its Applications. Cambridge University Press, Cambridge, expanded edition, 2014.
	
	\bibitem{Santalo}
	L. A. Santal$\acute{o}$. Integral geometry and geometric probability. Cambridge Mathematical Library. Cambridge
	University Press, Cambridge, second edition, 2004. With a foreword by Mark Kac.
	
	\bibitem{XYZZ2022}
	D. Xi, D. Yang, G. Zhang, and Y. Zhao. The Lp chord Minkowski problem. Advanced Nonlinear Studies, in
	press.
	
	\bibitem{GXZ2023}
	D. Xi, L. Guo, and Y. Zhao. The Lp chord Minkowski problem for $0\leqslant p<1.$ arXiv:2301.07603v1 [math.MG], in press.
	
	\bibitem{XL43}
	D. Xi and G. Leng, Dar’s conjecture and the log-Brunn-Minkowski inequality, J. Differential Geom. 103(2016), 145-189.
	
	\bibitem{H.Y45}
	H. Yagisita, Non-uniqueness of self-similar shrinking curves for an anisotropic curvature flow, Calc. Var. Partial Differential Equations 26(2006), 49-55.
	
	\bibitem{Zhu14}
	G. Zhu. The logarithmic Minkowski problem for polytopes. Adv. Math., 262:909–931, 2014.
		

\end{thebibliography}

\end{document}